\documentclass[11pt, reqno]{amsart}

\usepackage{graphicx}
\usepackage{epstopdf}
\usepackage{subcaption}
\usepackage[utf8]{inputenc}
\usepackage[T1]{fontenc}
\usepackage[english]{babel}
\usepackage{amsmath, amsthm}
\usepackage{ae}
\usepackage{icomma}
\usepackage{units}
\usepackage{color}
\usepackage{bbm}
\usepackage{caption}
\usepackage{array}
\usepackage[hmarginratio=1:1]{geometry}
\usepackage[hyphens]{url}
\usepackage[pdfpagelabels=false, hidelinks]{hyperref}

\usepackage{mathrsfs}
\usepackage{amssymb, amsfonts}
\usepackage{placeins} 
\usepackage{tikz}
\usetikzlibrary{patterns}
\usepackage{float}
\usepackage[nodayofweek]{datetime}
\usepackage{enumitem}
\usepackage{mathtools}
\usepackage[numbers, sort]{natbib}
\usepackage{dsfont}
\usepackage{comment}
\usepackage{csquotes}

\newcommand{\R}{\ensuremath{\mathbb{R}}}
\newcommand{\N}{\ensuremath{\mathbb{N}}}

\newcommand{\eps}{\ensuremath{\varepsilon}}

\DeclareMathOperator{\Tr}{Tr}

\newtheorem{theorem}{Theorem}[section]

\newtheorem{lemma}[theorem]{Lemma}
\newtheorem{proposition}[theorem]{Proposition}
\newtheorem{corollary}[theorem]{Corollary}

\numberwithin{theorem}{section}
\theoremstyle{remark}
\newtheorem{remark}[theorem]{Remark}

\numberwithin{equation}{section}

\begin{document}

\title[On the blow-up of the one-particle Lieb-Thirring constant in dimension two]{On the blow-up of the one-particle Lieb-Thirring constant in dimension two}

\author{Matthias Baur}
\address{\textnormal{(M. Baur)} Institute of Analysis, Dynamics and Modeling, Department of Mathematics, University of
Stuttgart, Pfaffenwaldring 57, 70569 Stuttgart, Germany.}
\email{\href{mailto:matthias.baur@mathematik.uni-stuttgart.de}{matthias.baur@mathematik.uni-stuttgart.de}}

\begin{abstract}
We establish three-term asymptotics for the one-particle Lieb-Thirring constant in dimension two as $\gamma\to 0$ and identify, up to rescaling and symmetries, the profile of the asymptotically worst potential. The proof uses the well-known relation between the one-particle constant and an integral moment of optimizers of the Gagliardo-Nirenberg inequality. Through the optimizers' characterization as positive, vanishing solutions of the nonlinear elliptic equation $-\Delta Q + Q = Q^{s}$ in $\mathbb{R}^2$, the essence of the proof lies in the asymptotic analysis of solutions $Q$ as $s\to\infty$.  
\end{abstract}

\maketitle


\section{Introduction and main result}

Let $V\ge0$, $V\in L^{\gamma+d/2}(\R^d)$, and let $\lambda_1\le\lambda_2\le\dots<0$ be the negative eigenvalues of the Schrödinger operator $-\Delta-V$ on $L^2(\R^d)$, counted with multiplicity. For $\gamma\geq 0$, define
\begin{align*}
   \Tr(-\Delta-V)_-^\gamma := \begin{cases}  
   \# \{ n : \lambda_n <0 \} & \text{if } \gamma = 0, \\
   \sum_j|\lambda_j|^{\gamma}  & \text{if } \gamma > 0,
   \end{cases}
\end{align*}
and let
\begin{equation} \label{eq:LTconstant}
   L_{\gamma,d}:=\sup_{V\neq 0}\frac{\Tr(-\Delta-V)_-^\gamma}{\int_{\R^d}V^{\gamma+d/2}\, dx}.
\end{equation}
The constant $L_{\gamma,d}$ is the sharp Lieb-Thirring constant for eigenvalue moments of order $\gamma$. It is the optimal constant such that
\begin{align*}
    \Tr(-\Delta-V)_-^\gamma  \leq  L_{\gamma,d} \int_{\R^d}V^{\gamma+d/2}\, dx 
\end{align*}
for any $V\in L^{\gamma+d/2}(\R^d)$. It is well-known that $L_{\gamma,d}$ is finite if and only if
\begin{align*}
    \gamma \geq 1/2 \quad &\text{ if } d=1, \\
    \gamma  > 0 \qquad &\text{ if } d = 2,\\
    \gamma \geq 0 \qquad  &\text{ if } d\geq 3,
\end{align*}
see \cite{LiebThirring1975,LiebThirring1976,Weidl1996,Cwikel1977,Lieb1976,Rozenblum1972,Frank2021,FrankLaptevWeidl2023,Lieb1980, Rozenblum1976,Hundertmark2007}. The sharp one-particle constant $L^{(1)}_{\gamma,d}$ is defined for $\gamma \geq 0$ by
\begin{equation} \label{eq:L1}
L^{(1)}_{\gamma,d}:=\sup_{V\neq 0}\frac{|\lambda_1|^{\gamma}}{\int_{\R^d}V^{\gamma+d/2}\, dx},
\end{equation}
and it is likewise the best constant such that
\begin{align} \label{eq:L1_bound}
    |\lambda_1|^{\gamma} \leq  L^{(1)}_{\gamma,d} \int_{\R^d}V^{\gamma+d/2}\, dx 
\end{align}
for any $V\in L^{\gamma+d/2}(\R^d)$, see Keller \cite{Keller1961}. By definition, it follows
\begin{align} \label{eq:LTconstant_geq_L1}
    L_{\gamma,d} \geq  L^{(1)}_{\gamma,d}
\end{align}
for any $\gamma> 0$ and $d\in\N$. Determining the value of the sharp constants $L_{\gamma,d}$ has been a decades-long endeavor, and while some cases have been covered by now and there have also been important recent developments \cite{Read2026,OpenAI2026,LaptevWeidl2000,HundertmarkLiebThomas1998}, the problem remains unsolved for $\gamma<3/2$ and $d \geq 2$. 

Here, we consider the case $d=2$. In two dimensions, the exact values of the sharp constants $L_{\gamma,2}$ are unknown for $\gamma \in (0,\frac{3}{2})$. It is known that $L_{\gamma,2} >  L^{(1)}_{\gamma,2}$ for $\gamma>1$, see \cite{FrankGontierLewin2021}, but it is unclear if $L_{\gamma,2} =  L^{(1)}_{\gamma,2}$ or $L_{\gamma,2} >  L^{(1)}_{\gamma,2}$ for $\gamma\in(0,1]$. The Lieb-Thirring conjecture predicts that $L_{\gamma,2} =  L^{(1)}_{\gamma,2}$ for $\gamma\in(0,1]$. At the left end point of the interval, it is known that $L_{0,2}$ and $L^{(1)}_{0,2}$ are infinite since arbitrarily weak potentials $V$ can cause at least one bound state \cite[Theorem XIII.11]{ReedSimon1978}, see also \cite{Simon1976}. In this paper, we investigate the behavior of the Lieb-Thirring constants, in particular the one-particle constant, in the limit $\gamma \to 0$. This is a regime that has received very little attention so far, cf.~\cite[Section 4.5]{Frank2021}. The main result of this article is the following.

\begin{theorem} \label{thm:L1thm}
   \begin{align*}
    L^{(1)}_{\gamma,2} = \frac{1}{4\pi e}\bigg(\frac{1}{\gamma}+\ln\frac{1}{\gamma}+(1-2\gamma_E) + o(1)\bigg), \qquad \gamma \to 0,
\end{align*} 
where $\gamma_E$ is the Euler-Mascheroni constant, and for the corresponding optimal potential $V_\gamma$ in \eqref{eq:L1} holds, up to rescaling and symmetries, locally uniformly
\begin{align*}
    V_\gamma(x) =  4 \gamma  (1+|x|^2)^{-2} (1+o(1)),  \qquad \gamma\to 0.
\end{align*}
\end{theorem}

The proof of Theorem \ref{thm:L1thm} uses the well-known duality of the one-particle inequality to the Gagliardo-Nirenberg inequality. The optimal constant $L^{(1)}_{\gamma,2}$ is attained by $V=Q^{s-1}$ where $Q$ is the optimizer of the Gagliardo-Nirenberg inequality. The optimizer $Q$ is characterized by the Euler-Lagrange equation $-\Delta Q + Q = Q^{s}$ in $\R^2$ with $s=1+2/\gamma$. Since $s\to \infty$ as $\gamma \to 0$, one then needs to compute asymptotics for $Q$ as $s\to \infty$. The analysis of $Q$ for large $s$ is inspired by the work of Ren and Wei \cite{RenWei1994} on least-energy solutions of $-\Delta u=u^p$ as $p\to\infty$ and by the blow-up analysis of $-\Delta u + \lambda u = u^{p}$  on bounded planar domains of Adimurthi and Grossi \cite{AdimurthiGrossi2004} and El Mehdi and Grossi \cite{ElMehdiGrossi2004}, in particular the identification of the limit profile and proof that $\Vert u \Vert_{\infty} \to \sqrt{e}$ as $p \to \infty$. One notable difference from bounded domains is that in the setting of the present paper, the optimizer $Q$ on $\R^2$ is radial and one can directly exploit comparison estimates on the level of the radial ODE corresponding to the Euler-Lagrange equation. 

Regarding the true Lieb-Thirring constant $L_{\gamma,2}$, we mention that Lieb's method \cite{Lieb1976,Lieb1980} yields, for $V\geq 0$, $E\geq 0$ and any convex, polynomially bounded function $f\geq 0$ with $f(0)=0$,
\begin{align*}
    \# \{ n : \lambda_n < -E \}\leq \frac{1}{c_f} \int_0^\infty \frac{1}{t}\, e^{-tE} \int_{\R^2} \frac{f(tV(x))}{4\pi t}\, dx \, dt,
    \qquad c_f = \int_0^\infty f(s)\, s^{-1} e^{-s}\, ds .
\end{align*}
Choosing $f(s)=(s-a)_+$, $a>0$, for which $c_f = e^{-a}-aE_1(a)$ with $E_1(a)=\int_a^\infty s^{-1}e^{-s}\,ds$, and integrating via $\Tr(-\Delta-V)_-^\gamma = \gamma\int_0^\infty E^{\gamma-1} \# \{ n : \lambda_n < -E \} \, dE$, one obtains
\[
L_{\gamma,2}\le  \frac{\Gamma(1+\gamma)\,a^{-\gamma}}{4\pi\gamma(1+\gamma) (e^{-a}-aE_1(a) )}, \qquad a>0 .
\]
Setting $a=\gamma$, we have $a^{-\gamma} \to 1$ and $(e^{-a}-aE_1(a) ) \to 1$ as $\gamma\to 0$, so 
$$L_{\gamma,2}\le \frac{1}{4\pi\gamma}(1+o(1)), \qquad \gamma\to 0.$$
We conclude with Theorem \ref{thm:L1thm} that $L_{\gamma,2}$ blows up like $\gamma^{-1}$, but we are unable to determine the exact constant in front of $\gamma^{-1}$.

\begin{corollary}
\begin{align*}
   \frac{1}{4\pi e  \gamma} (1+o(1)) \leq  L_{\gamma,2} \leq  \frac{1}{4\pi  \gamma} (1+o(1)), \qquad \gamma \to 0.
\end{align*}
\end{corollary}

Closing the gap of a factor $e$ remains an open problem.

\section{A priori lower bound on $L^{(1)}_{\gamma,2}$}

Instead of diving straight into the proof of the main theorem, we start with a simple proof of the following bound on the one-particle constant $L^{(1)}_{\gamma,2}$. 

\begin{proposition} \label{prop:L1-lowerbnd}
    \begin{align*}
    L^{(1)}_{\gamma,2} \ge \frac{(1+\gamma)^{1+\gamma}(1+2\gamma)}{4\pi \gamma^{1+\gamma}}\,e^{-1-2\gamma(1+\gamma_E)}, \qquad \gamma>0.
\end{align*}
\end{proposition}

To obtain this bound, it is enough to use a particular trial potential. The proof is instructive, since it motivates what the asymptotically worst potential looks like, see Remark \ref{rem:motivation}, and by elementary estimates, one sees that the proposition already contains one side of the asymptotic expansion of $L^{(1)}_{\gamma,2}$.

\begin{corollary}
   \begin{align*}
    L^{(1)}_{\gamma,2} \ge \frac{1}{4\pi e}\bigg(\frac{1}{\gamma}+\ln\frac{1}{\gamma}+(1-2\gamma_E)-2\gamma\bigg), \qquad \gamma>0.
\end{align*} 
\end{corollary}

\begin{proof}[Proof of Proposition \ref{prop:L1-lowerbnd}]
Let $V= \mu \varphi$ with $\varphi \geq 0$, $\varphi \in L^1(\R^2) \cap L^\infty(\R^2)$, $\int_{\R^2 } \varphi \, dx = 1$ and $\mu>0$. By the Birman-Schwinger principle, we have $\lambda_1(-\Delta - V) \leq -\kappa^2$ if and only if $\lambda_{\max}( B_\kappa) \geq 1$ where
\begin{align*}
    B_\kappa = V^{1/2} (-\Delta + \kappa^2)^{-1} V^{1/2}, \qquad \kappa>0
\end{align*}
denotes the Birman-Schwinger operator. It is compact since $\varphi \in L^1(\R^2) \cap L^\infty(\R^2)$ and in fact even Hilbert-Schmidt. Moreover, its integral kernel is
\begin{align*}
    B_\kappa(x,y) = \frac{1}{2\pi} V^{1/2}(x) K_0(\kappa |x-y|) V^{1/2}(y), \qquad x,y\in \R^2,
\end{align*}
where $K_0$ is a modified Bessel function of the second kind. Taking the trial function $\varphi^{1/2} \in L^2(\R^2)$ and then using the estimate 
\begin{align*}
    K_0(\rho)\;\ge\;\ln\frac{2}{\rho}-\gamma_E\qquad\text{for all } \rho>0 ,
\end{align*}
where $\gamma_E$ is the Euler-Mascheroni constant, we obtain
\begin{align*}
    \lambda_{\max}(B_\kappa)\ &\ge\ \frac{\mu}{2\pi}\iint_{\R^2 \times \R^2} \varphi(x)K_0(\kappa|x-y|)\varphi(y)\,dx\,dy\ \\
    &\ge\ \frac{\mu}{2\pi}\Big(\ln\frac{2}{\kappa}-\gamma_E+I(\varphi)\Big),
\end{align*}
where 
\begin{align*}
 I(\varphi)=\iint_{\R^2 \times \R^2}\varphi(x)\varphi(y)\ln\frac{1}{|x-y|} \, dx \, dy .
\end{align*}
Therefore, the Birman-Schwinger principle implies 
\begin{align*}
    |\lambda_1(-\Delta-V)|\ \ge\ 4 e^{-2\gamma_E + 2I(\varphi) -4\pi/\mu}.
\end{align*}
For the one-particle constant, this leads to the lower bound
\begin{align} \label{eq:chain-oneparticle-est}
\begin{split}
       L^{(1)}_{\gamma,2}\ &=\ \sup_{ V\geq 0} \frac{|\lambda_1(-\Delta-V)|^\gamma}{\int_{\R^2} V^{1+ \gamma} \, dx}\\
    &\geq \sup_{\mu >0 } \frac{|\lambda_1(-\Delta-\mu\varphi)|^\gamma}{\mu^{1+\gamma}\int_{\R^2} \varphi^{1+ \gamma} \, dx} \\
    &\geq \sup_{\mu >0 } \frac{(4 e^{-2\gamma_E + 2I(\varphi) -4\pi/\mu})^\gamma}{\mu^{1+\gamma}\int_{\R^2} \varphi^{1+ \gamma} \, dx} \\
    & =  \frac{e^{-(1+\gamma)}(1+\gamma)^{1+\gamma} \big(4e^{-2\gamma_E}\big)^{\gamma}}{(4\pi\gamma)^{1+\gamma}}\cdot\frac{e^{2\gamma I(\varphi)}}{\int_{\R^2 }\varphi^{1+\gamma} \, dx}  
\end{split}
\end{align}
where to obtain the last line we have optimized with respect to $\mu>0$. The supremum is attained at $\mu_*=\frac{4\pi \gamma}{(1+\gamma)}$. It remains to determine a good potential profile $\varphi$. We choose 
$$\varphi_0(x)= \frac{1}{\pi}(1+|x|^2)^{-2}.$$ 
Direct evaluation of the integrals then gives
\begin{align*}
    I(\varphi_0)=-\frac12\int_1^\infty\frac{\ln u}{u^2}\,du=-\frac12,\qquad \int_{\R^2 } \varphi_0^{1+\gamma}\, dx=\frac{\pi^{-\gamma}}{1+2\gamma}
\end{align*}
and therefore
\begin{align*}
    \frac{e^{2\gamma I(\varphi_0)}}{\int_{\R^2 }\varphi_0^{1+\gamma} \, dx}  = (1+2\gamma) e^{\gamma (-1+ \ln \pi)}.
\end{align*}
Inserting this into \eqref{eq:chain-oneparticle-est} proves the bound of the proposition. 
\end{proof}

\begin{remark}
In the argument above, the trial potential $V$ chosen to estimate $L^{(1)}_{\gamma,2}$ from below is precisely
\begin{align*}
    V(x) = \mu_* \varphi_0(x) = \frac{4\gamma }{(1+\gamma)}  (1+|x|^2)^{-2}  ,
\end{align*}
which Theorem \ref{thm:L1thm} predicts to be asymptotically optimal. 
\end{remark}

\begin{remark} \label{rem:motivation}
We want to briefly motivate the choice of $\varphi_0$ from the argument in the proof above. Taking the logarithm of the last factor from the last line of \eqref{eq:chain-oneparticle-est}, we have
\begin{align*}
    \ln\bigg(  \frac{e^{2\gamma I(\varphi)}}{\int_{\R^2 }\varphi^{1+\gamma} \, dx} \bigg) = 2 \gamma I(\varphi) - \ln \bigg(\int_{\R^2 }\varphi^{1+\gamma} \, dx \bigg).
\end{align*}
By Jensen's inequality, 
\begin{align*}
   - \ln \bigg( \int_{\R^2 }\varphi^{1+\gamma} \, dx \bigg) \leq -\gamma \int_{\R^2 }\varphi \ln \varphi \, dx,
\end{align*}
thus
\begin{align*}
    \ln\bigg(  \frac{e^{2\gamma I(\varphi)}}{\int_{\R^2 }\varphi^{1+\gamma} \, dx} \bigg) \leq  2\gamma \bigg(  I(\varphi) - \frac{1}{2}\int_{\R^2 }\varphi \ln \varphi \, dx \bigg) .
\end{align*}
The logarithmic Hardy-Littlewood-Sobolev inequality by Carlen and Loss \cite{Carlen1992} in dimension two shows that 
\begin{align} \label{eq:logHLS}
    I(\varphi) \leq \frac{1}{2}\int_{\R^2 }\varphi \ln \varphi \, dx + \frac{1}{2} (1+ \ln \pi)
\end{align}
for any profile $\varphi \in L^1(\R^2) \cap L^\infty(\R^2)$ with $\int_{\R^2} \varphi \, dx =1$, $\int_{\R^2} \varphi  \ln(1+|x|^2)\, dx <\infty$, so it follows the upper bound
\begin{align} \label{eq:quotient_optimalbound}
    \frac{e^{2\gamma I(\varphi)}}{\int_{\R^2 }\varphi^{1+\gamma} \, dx} \leq e^{\gamma(1+\ln \pi)}
\end{align}
for any $\varphi \in L^1(\R^2) \cap L^\infty(\R^2)$, $\varphi \geq 0$, $\int_{\R^2} \varphi \, dx = 1$\footnote{Note when $\int_{\R^2} \varphi  \ln(1+|x|^2)\, dx =\infty$, then $I(\varphi)=-\infty$, so the left-hand side is zero and the bound is trivial.}. The best we can hope for with any potential profile $\varphi$ is thus equality in \eqref{eq:quotient_optimalbound}. Equality in the logarithmic Hardy-Littlewood-Sobolev inequality \eqref{eq:logHLS} is attained by $\varphi_0$. This is what motivates the choice $\varphi = \varphi_0$ in the proof. We remark that while $\varphi_0$ is not optimal for \eqref{eq:quotient_optimalbound}, it is almost optimal in the sense that it achieves the upper bound \eqref{eq:quotient_optimalbound} up to a loss of $O(\gamma^2)$ as $\gamma\to 0$. 
\end{remark}

\begin{remark}
Theorem \ref{thm:L1thm} is in fact stronger than the logarithmic Hardy-Littlewood-Sobolev inequality \eqref{eq:logHLS}. Indeed, for any profile $\varphi \in L^1(\R^2) \cap L^\infty(\R^2)$ with $\int_{\R^2} \varphi \, dx =1$, $\int_{\R^2} \varphi  \ln(1+|x|^2)\, dx <\infty$ follows by \eqref{eq:chain-oneparticle-est}, Theorem \ref{thm:L1thm} and the estimate $\frac{e^{\gamma}}{(1+\gamma)^{1+\gamma}  } \leq 1$ for $\gamma\geq 0$,
\begin{align*}
    e^{2\gamma I(\varphi)} &\leq (4\pi e \gamma L_{\gamma,2}^{(1)} ) (\pi\gamma)^\gamma  e^{2\gamma_E \gamma}\frac{e^{\gamma}}{(1+\gamma)^{1+\gamma}  } {\int_{\R^2 }\varphi^{1+\gamma} \, dx}  \\
    &\leq (1+\gamma\ln\frac{1}{\gamma}+(1-2\gamma_E)\gamma + o(\gamma))  e^{\ln \pi -\gamma\ln\frac{1}{\gamma}+2\gamma_E \gamma} {\int_{\R^2 }\varphi^{1+\gamma} \, dx}.
\end{align*}
Taking the logarithm of both sides and letting $\gamma\to 0$ gives
\begin{align*}
    2I(\varphi) \leq (1+ \ln \pi) + \lim_{\gamma\to 0} \frac{1}{\gamma} \ln \bigg( {\int_{\R^2 }\varphi^{1+\gamma} \, dx} \bigg)
\end{align*}
and by dominated convergence,
\begin{align*}
    \lim_{\gamma\to 0} \frac{1}{\gamma} \ln \bigg( {\int_{\R^2 }\varphi^{1+\gamma} \, dx} \bigg) =\lim_{\gamma\to 0} \frac{1}{\gamma} \bigg( {\int_{\R^2 }\varphi^{1+\gamma} \, dx} -1\bigg) &= \lim_{\gamma\to 0}  {\int_{\R^2 }\varphi \,\frac{\varphi^\gamma-1}{\gamma} \, dx} \\
    &= {\int_{\R^2 }\varphi \ln \varphi \, dx},
\end{align*}
which shows \eqref{eq:logHLS}.

\end{remark}

\section{Asymptotic expansion of $L^{(1)}_{\gamma,2}$}

In this section, we prove Theorem \ref{thm:L1thm}.

\subsection{The Gagliardo-Nirenberg optimizer}

For $\gamma>0$ let
\begin{equation*}
  p :=  1+\frac1\gamma,\qquad s := 2p-1 = 1+\frac{2}{\gamma}.
\end{equation*}
We recall that the one-particle constant is the best constant in \eqref{eq:L1_bound}. It is dual to the
Gagliardo-Nirenberg inequality
\begin{equation}\label{eq:GN}
  \int_{\R^2}|\nabla u|^2\,dx\ \ge\ K^{(1)}_{p,2}\Big(\int_{\R^2}|u|^{2p}\,dx\Big)^{\frac1{p-1}}\|u\|_{L^2(\R^2)}^{-\frac{2}{p-1}},
  \qquad u\in H^1(\R^2),
\end{equation}
in the sense that the optimal constants satisfy \cite[Sect.~4.2]{Frank2021}, \cite{LiebThirring1976}
\begin{equation}\label{eq:duality}
  L^{(1)}_{\gamma,2}\,K^{(1)}_{p,2} = \frac{1}{\gamma+1}\Big(\frac{\gamma}{\gamma+1}\Big)^{\gamma}.
\end{equation}

\begin{proposition} \label{prop:Q}
Let $\gamma>0$.
\begin{enumerate}
\item The Gagliardo-Nirenberg inequality \eqref{eq:GN} has an optimizer, which is radial and unique up to translation,   dilation and multiplication by a constant. After these symmetries the optimizer $Q$ is
  positive, centered at the origin and solves the Euler-Lagrange equation
  \begin{equation}\label{eq:EL}
    -\Delta Q + Q = Q^{s}\qquad\text{in } \R^2 ,
  \end{equation}
  and it is the unique positive solution of \eqref{eq:EL} that vanishes at infinity. Moreover $Q$ is smooth, radially
  strictly decreasing, $Q(0)=\|Q\|_{L^\infty(\R^2)}$, and $Q(0)>1$.
  \item There are $0<c_\gamma<C_\gamma$ with
  \begin{equation}\label{eq:decay}
    c_\gamma\,(1+|x|)^{-1/2}e^{-|x|}\le Q(x)\le C_\gamma\,(1+|x|)^{-1/2}e^{-|x|}
  \end{equation}
  and
  \begin{align}
    Q(x)&=\frac{1}{2\pi}\int_{\R^2}K_0(|x-y|)\,Q(y)^{s}\,dy \label{eq:green} 
  \end{align}
  where $K_0$ is the modified Bessel function of the second kind.
\item The optimizer $Q$ satisfies
  \begin{align}
    \int_{\R^2} Q^{s+1}\, dx &=  \frac{s+1}{2}\, \int_{\R^2} Q^{2} \, dx , \label{eq:pohozaev} \\
    \int_{\R^2} Q^{s}\, dx&= \int_{\R^2} Q \, dx,  \label{eq:pohozaev2} \\
    \int_{\R^2} |\nabla Q|^{2} \, dx &= \frac{1}{\gamma}\, \int_{\R^2} Q^{2} \, dx  .  \label{eq:pohozaev3}
  \end{align}
\item It holds
  \begin{equation}\label{eq:L1Q}
    L^{(1)}_{\gamma,2} = \bigg( \int_{\R^2} Q^{s+1}\, dx \bigg)^{-1}
  \end{equation}
  and the supremum in \eqref{eq:L1} is attained, uniquely up to rescaling and symmetries, for
  $V=Q^{s-1}$. 
\end{enumerate}
\end{proposition}

\begin{proof}[References and proof sketch]
(1) We refer to \cite{Frank2021} and \cite{FrankGontierLewin2021} and the references therein. Detailed proofs are
collected in \cite{Frank2013notes}. The fact that $Q(0)>1$ follows from $\int_{\R^2} |\nabla Q|^2 \, dx = \int_{\R^2} Q^2 (Q^{s-1}-1) \, dx >0$ which requires that $\Vert Q \Vert_{L^{\infty}(\R^2)}>1$.

(2) The two-sided bound \eqref{eq:decay} is a consequence of \cite[Proposition 4.1]{GidasNiNirenberg1981}. It is used in this form in \cite[Sect.~6.2]{Frank2021}. Since $Q\in H^1(\R^2)$ and $Q^{s}\in L^1(\R^2)$, equation
\eqref{eq:EL} is equivalent to $Q=(-\Delta+1)^{-1}Q^{s}$ and the integral kernel of $(-\Delta+1)^{-1}$ on $\R^2$ is $\frac{1}{2\pi}K_0(|x-y|)$.

(3) For \eqref{eq:pohozaev} and \eqref{eq:pohozaev3} we refer to \cite[(35)--(36)]{Frank2021}. Integrating   \eqref{eq:green} and Fubini together with $\int_{\R^2} K_0(|x|) \, dx = 2\pi$ shows \eqref{eq:pohozaev2}.

(4) Let $m=\int Q^2 dx$. Evaluating \eqref{eq:GN} at $Q$ and using \eqref{eq:pohozaev} and \eqref{eq:pohozaev3},
\[
  K^{(1)}_{p,2}=\frac{(p-1)m\cdot m^{\frac1{p-1}}}{(pm)^{\frac1{p-1}}}=(p-1)\,p^{-\frac1{p-1}}\,m ,
\]
while the right-hand side of \eqref{eq:duality} equals $\frac{p-1}{p} p^{-\frac1{p-1}}$. Dividing gives 
$$L^{(1)}_{\gamma,2}=\frac{1}{pm}=\bigg(\int_{\R^2} Q^{s+1} dx \bigg)^{-1},$$
cf.\ \cite[Sect.~6.2]{Frank2021}. The duality argument also shows that optimizers of
\eqref{eq:GN} and optimal potentials in \eqref{eq:L1} are in one-to-one correspondence
modulo symmetries \cite[Sect.~4.2]{Frank2021}. If $V=Q^{s-1}$, then $Q$ is the ground state of $-\Delta-V$
with eigenvalue $-1$ by \eqref{eq:EL}, and $V^{\gamma+1}=Q^{s+1}$. 
\end{proof}

In view of \eqref{eq:green}, let us record a few properties of the modified Bessel function $K_0$ which will be of use later on.

\begin{lemma} \label{lem:K0}
The function $K_0$ is positive and strictly decreasing on $(0,\infty)$, $\int_{\R^2}K_0(|x|)\, dx = 2\pi$, and there exists a constant $C$ such that for $0<\rho\le 1$
\begin{align} \label{eq:K0-expansion}
    0 \le K_0(\rho) - \ln\frac{2}{\rho} + \gamma_E \le C\rho^2 \bigg(1+\ln\frac{1}{\rho} \bigg) .
\end{align}
Consequently, $K_0(\rho) \le \ln^+\frac{2}{\rho} + C$ for all $\rho>0$, where $\ln^+ = \max\{\ln, 0\}$.
\end{lemma}

\begin{proof}
Positivity and monotonicity of $K_0$ follow directly from the integral representation 
$$K_0(\rho)=\int_0^\infty e^{-\rho\cosh t}\,dt,$$
see \cite[(10.32.9)]{DLMF}, and the value of the integral is standard. By the series representation \cite[(10.31.2)]{DLMF}
\begin{align*}
    K_0(\rho) = -\Big(\ln\frac{\rho}{2} + \gamma_E\Big) I_0(\rho) + \sum_{k=1}^\infty \Big(\frac{\rho^2}{4}\Big)^k \frac{H_k}{(k!)^2}, \qquad H_k = \sum_{j=1}^k \frac1j,
\end{align*}
where $I_0(\rho) = \sum_{k\ge0} (\rho^2/4)^k/(k!)^2$ is the modified Bessel function of the first kind, we have
\begin{align*}
    K_0(\rho) - \ln\frac{2}{\rho} + \gamma_E = \Big(\ln\frac{2}{\rho} - \gamma_E\Big)\big(I_0(\rho)-1\big) + \sum_{k=1}^\infty \Big(\frac{\rho^2}{4}\Big)^k \frac{H_k}{(k!)^2} .
\end{align*}
For $0<\rho\le1$ both terms are nonnegative, since $\ln 2 - \gamma_E >0$, and $0\le I_0(\rho)-1\le C\rho^2$, which gives \eqref{eq:K0-expansion}. For $\rho\ge1$ we use $K_0(\rho)\le K_0(1) < 1$.
\end{proof}

\subsection{Rescaled picture}

The main task is to find asymptotics for the integral
\begin{align*}
     \int_{\R^2} Q^{s+1} \, dx
\end{align*}
of the Gagliardo-Nirenberg optimizer $Q$ as $s \to \infty$, since this determines the asymptotics of the one-particle constant $L_{\gamma,2}^{(1)}$ by \eqref{eq:L1Q}. To obtain these, we start by rescaling around the maximum of $Q$ at $x=0$. Set
\begin{align} \label{eq:def_rescaled_var}
    \varepsilon^{2}:=\frac{1}{s\,Q(0)^{s-1}},\qquad U_s(z):=s\bigg(\frac{Q(\varepsilon z)}{Q(0)}-1\bigg). 
\end{align}
Then,
\begin{align*}
    Q(\varepsilon z) &= Q(0)\bigg(1+\frac{U_s(z)}{s}\bigg), \qquad 
    Q(\varepsilon z)^{s}=\frac{Q(0)}{s\varepsilon^{2}}\bigg(1+\frac{U_s(z)}{s}\bigg)^{s},
\end{align*}
and $U_s$ is radial, $U_s(0)=0$, $\nabla U_s(0)=0$, $-s<U_s\le 0$, and satisfies
\begin{align} \label{eq:Qs-rescaled}
     0 \le \bigg(1+\frac{U_s(z)}{s}\bigg)^{s} \le e^{U_s(z)} \le 1 .
\end{align}
The Euler-Lagrange equation $-\Delta Q + Q = Q^s$ gets translated to the equation
\begin{equation}\label{eq:innerEL}
-\Delta U_s=\bigg(1+\frac{U_s}{s}\bigg)^{s}- s\varepsilon^2 \bigg(1+\frac{U_s}{s}\bigg)\quad\text{in }\R^2 .
\end{equation}
Here, the mass term $+Q$ of $-\Delta Q + Q = Q^s$ became the forcing $s\eps^2(1+U_s/s)$ on the right-hand side. It will turn out that $s\eps^2 = Q(0)^{1-s}$ and thus the forcing term are exponentially small, since $Q(0)$ converges to a constant strictly above one as $s\to\infty$. The elementary limit $(1+\frac{u}{s})^{s} \to e^u$ for any $u \in \R$ motivates the idea that $U_s$ might converge in some sense to a function $U$ that solves the limiting equation $-\Delta U = e^U$. So let us introduce $U:\R^2 \to \R$ defined by 
    \begin{align*}
    U(z)=-2\ln\Big(1+\frac{|z|^{2}}{8}\Big).
\end{align*}
Then $U$ satisfies $-\Delta U=e^{U}$, $U(0)=0$, $\nabla U (0)=0$. It will play in the following the role of a limiting profile for $U_s$.

Let us now write
\begin{align} \label{eq:Qsplus1_IJ}
    \int_{\R^2} Q^{s+1} \, dx = \frac{Q(0)^2}{s} \int_{\R^2} \Big(1+\frac{U_s}{s}\Big)^{s+1} \, d z = \frac{Q(0)^2}{s} \bigg( I(s) + \frac{1}{s} J(s)  \bigg)
\end{align}
where

\begin{align*}
I(s)&:=\int_{\R^2} \Big(1+\frac{U_s}{s}\Big)^{s} \, d z, \\
J(s)&:=\int_{\R^2} U_s\Big(1+\frac{U_s}{s}\Big)^{s} \, d z 
\end{align*}
and let us further denote
\begin{align*}
    L(s)&:=\int_{\R^2}\ln|z|\, \Big(1+\frac{U_s}{s}\Big)^{s}\, d z.
\end{align*}

We claim that the integrals $I$, $J$ and $L$ converge to respective integrals of the limiting profile $U$.
\begin{proposition} \label{prop:IJL_limits}
As $s\to\infty$,
\begin{align}
    I(s) &\to  \int_{\mathbb{R}^2}e^{U} \, dz = 8 \pi \label{eq:I_limit}, \\
    J(s) &\to  \int_{\mathbb{R}^2} Ue^{U} \, dz = -16\pi \label{eq:J_limit}, \\
    L(s) &\to  \int_{\mathbb{R}^2} \ln|z| \; e^{U} \, dz =12\pi\ln 2 \label{eq:L_limit}.
\end{align}
\end{proposition}

We postpone the proof of this proposition to a later section. To prove Theorem \ref{thm:L1thm} at leading order it would be enough to have \eqref{eq:I_limit}, \eqref{eq:J_limit} and asymptotics for $Q(0)$, since equation \eqref{eq:Qsplus1_IJ} directly connects asymptotics for the integrals $I$, $J$ with those of $\int_{\R^2} Q^{s+1} \, dx$ and thus $L_{\gamma,2}^{(1)}$. However, in order to go beyond leading order, we can later exploit two relations between the integrals $I$, $J$ and $L$.

\subsection{Two identities between integrals}

Two identities relate the integrals $I$, $J$ and $L$ from Proposition \ref{prop:IJL_limits}. The first one is a consequence of the Green representation \eqref{eq:green} evaluated at the origin.
 
\begin{lemma} \label{lem:I}
For every $s$,
\begin{align} \label{eq:I}
    2\pi s = \bigg(\ln\frac{2}{\eps}-\gamma_E\bigg) I(s) - L(s) + R_1(s),
\end{align}
where
\begin{align} \label{eq:R1-def}
    R_1(s) := \int_{\R^2} \bigg[K_0(\eps|z|) - \ln\frac{2}{\eps|z|}+\gamma_E\bigg] \bigg(1+\frac{U_s(z)}{s}\bigg)^{s} \, dz .
\end{align}
\end{lemma}
 
\begin{proof}
Equation \eqref{eq:green} at $x=0$ with the substitution $y=\eps z$ gives
\begin{align*}
    2\pi Q(0) = \int_{\R^2}K_0(|y|)Q(y)^s\,dy = \eps^2 \,\frac{Q(0)}{s\eps^2} \int_{\R^2} K_0(\eps|z|) \bigg(1+\frac{U_s(z)}{s}\bigg)^{s} \, dz 
\end{align*}
and therefore
\begin{align} \label{eq:I-exact}
    2\pi s = \int_{\R^2} K_0(\eps|z|) \bigg(1+\frac{U_s(z)}{s}\bigg)^{s} \, dz .
\end{align}
The lemma follows now from writing $K_0(\eps|z|) = \ln\frac2\eps - \ln|z| - \gamma_E + [...]$.
\end{proof}
 
The second identity follows from \eqref{eq:pohozaev} and captures the fact that $\int_{\R^2} Q^2 \, dx$ is asymptotically carried by a Bessel tail.
 
\begin{lemma} \label{lem:P}
For every $s$,
\begin{align} \label{eq:P-inner}
    I(s) + \frac{1}{s}J(s) = \frac{s+1}{8\pi s}\, I(s)^2 \,  (1+\delta(s) )^2
\end{align}
and
\begin{align}
    \int_{\R^2} Q^{2} \, dx &= \frac{1}{4\pi} \bigg( \int_{\R^2} Q^{s} \, dx \bigg)^2\,(1+\delta(s) )^2, \label{eq:P-L2}\\
    \int_{\R^2} Q^{s+1} \, dx &= \frac{s+1}{8\pi}\,\bigg( \int_{\R^2} Q^{s} \, dx \bigg)^2 (1+\delta(s) )^2 ,\label{eq:P}
\end{align}
where $|\delta(s)| \leq R_2(s) $ with
\begin{align} \label{eq:Ds-def}
    R_2(s) := \frac{2\sqrt{\pi}}{\int_{\R^2} Q^{s} \, dx} \,\bigg\Vert Q - \frac{\int_{\R^2} Q^{s} \, dx}{2\pi} K_0(|\cdot|) \bigg\Vert_{L^2(\R^2)} .
\end{align}
\end{lemma}
 
\begin{proof}
Reverse triangle inequality applied to
\begin{align*}
    Q(x) = \frac{\int_{\R^2} Q^{s} \, dx}{2\pi}K_0(|x|) + \bigg( Q(x) - \frac{\int_{\R^2} Q^{s} \, dx}{2\pi}K_0(|x|) \bigg) 
\end{align*}
gives
\begin{align*}
    \bigg| \Vert Q \Vert_{L^2(\R^2)} - \frac{\int_{\R^2} Q^{s} \, dx}{2\pi} \Vert K_0(|\cdot|)\Vert_{L^2(\R^2)} \bigg|  &\leq \bigg\Vert Q(x) - \frac{\int_{\R^2} Q^{s} \, dx}{2\pi}K_0(|x|)  \bigg\Vert_{L^2(\R^2)} \\
    &= \frac{\int_{\R^2} Q^{s} \, dx}{2\sqrt{\pi}}  R_2(s) .
\end{align*}
Since $\Vert K_0(|\cdot|)\Vert_{L^2(\R^2)}^2 = \pi$, 
\begin{align*}
    \Vert Q \Vert_{L^2(\R^2)} = \frac{1}{2\sqrt{\pi}} \int_{\R^2} Q^{s} \, dx \,  ( 1+ \delta(s))
\end{align*}
with $|\delta(s)| \leq R_2(s)$. Squaring yields \eqref{eq:P-L2}. Equation \eqref{eq:P} follows immediately from \eqref{eq:P-L2} and \eqref{eq:pohozaev}. Equation \eqref{eq:P-inner} is \eqref{eq:P} rewritten with the identities $\int_{\R^2}  Q^{s} dx = \frac{Q(0)}{s} I(s)$ and $\int_{\R^2} Q^{s+1} \, dx  = \frac{Q(0)^2}{s} ( I(s) + \frac{1}{s} J(s))$.
\end{proof}

Finally, we claim that the two remainders $R_1$ and $R_2$ are exponentially small as $s\to \infty$.
 
\begin{lemma}  \label{lem:remainders}
There is $c>0$ such that
\begin{align*}
    R_1(s) = O(e^{-cs}), \qquad R_2(s) = O(e^{-cs}) , \qquad s \to \infty.
\end{align*}
\end{lemma}

The exponential smallness of the remainders is what allows us to treat both identities in the asymptotic regime as if there were no remainder and gain subleading order terms for $I$, $\varepsilon$ and eventually $L_{\gamma,2}^{(1)}$.

\subsection{Proof of Proposition \ref{prop:IJL_limits} and Lemma \ref{lem:remainders}}
\label{subsec:proof-IJL}

In this subsection we prove the convergence of the integrals $I$, $J$, $L$ and the exponential smallness of the remainders $R_1$, $R_2$. The proof is organized as a chain of lemmas, each depending only on the previous ones. The starting point is a collection of crude bounds on $Q$ which allow us to show that $U$ is indeed the limiting profile of $U_s$. We then show two-sided bounds on $Q(0)$ and finally a decay estimate for $(1+U_s/s)^s$ which is uniform in $s$ on a ball of exponentially large radius and is obtained by comparison in the radial ODE. Proposition \ref{prop:IJL_limits} then follows by dominated convergence, and Lemma \ref{lem:remainders} by using said decay estimate.

Throughout this subsection, $C$, $C'$ and sometimes $C_1$,  $C_2$, ... usually denote positive constants independent of $s$ which may change from line to line, and $s$ is assumed to be large. We begin with one elementary lemma recording the exact values for the integrals of the limiting profile $U$.

\begin{lemma} \label{lem:U_integrals}
The function $U(z)=-2\ln(1+|z|^2/8)$ satisfies 
\begin{align*}
\int_{\mathbb{R}^2}e^{U} \, dz& =8\pi,\qquad&  \int_{\{|z| \leq r \}}e^{U} \, dz &= \frac{8\pi r^2 }{8+r^2}, \\
\int_{\mathbb{R}^2} Ue^{U} \, dz&=-16\pi,\qquad& \int_{\mathbb{R}^2} \ln|z|\,e^{U} \, dz&=12\pi\ln 2 .
\end{align*}
\end{lemma}

\begin{proof}
With the substitution $u=1+t^2/8$, $t\,dt = 4\,du$,
\begin{align*}
    \int_{\{|z|\le r\}} e^{U}\, dz = 2\pi \int_0^r \frac{t\, dt}{(1+t^2/8)^2} = 8\pi \int_1^{1+r^2/8} \frac{du}{u^2} = 8\pi \Big( 1 - \frac{8}{8+r^2} \Big) = \frac{8\pi r^2}{8+r^2},
\end{align*}
which gives the first two identities, and
\begin{align*}
    \int_{\R^2} Ue^{U}\, dz = -16\pi \int_1^\infty \frac{\ln u}{u^2}\, du = -16\pi .
\end{align*}
For the last identity substitute $t^2 = 8v$, so that $t\,dt = 4\,dv$ and $\ln t = \frac12 \ln(8v)$:
\begin{align*}
    \int_{\R^2} \ln|z|\, e^{U}\, dz = 4\pi \int_0^\infty \frac{\ln (8 v)}{(1+v)^2}\, dv = 4\pi \ln 8 + 4\pi \int_0^\infty \frac{\ln v}{(1+v)^2}\, dv = 12\pi \ln 2,
\end{align*}
since the last integral vanishes by the substitution $v\mapsto 1/v$.
\end{proof}

The analysis of $U_s$ rests on the following a priori bounds, which follow from the a priori lower bound on $L^{(1)}_{\gamma,2}$ established in Section 2.

\begin{lemma} \label{lem:crude}
There exist constants $C$ and $s_0$ such that for $s\ge s_0$,
\begin{enumerate}
\item $ \displaystyle   \int_{\R^2} Q^{s+1}\, dx  \le \frac{C}{s}$,
\item $\displaystyle  \int_{\R^2} Q^{2}\, dx \le \frac{C}{s^2}$,
\item $\displaystyle \int_{\R^2} |\nabla Q|^{2}\, dx \le \frac{C}{s}$, 
\item $Q(x) \le \dfrac{C}{s|x|}$ for all $x\ne0$,
\item $\displaystyle\int_{\R^2} Q^{s}\, dx \le \dfrac{C}{s}$,
\item $\displaystyle \int_{\R^2} (1+|x|+|\ln|x|| )\,Q(x)\,dx \le C$,
\item $I(s)  \le C$.
\end{enumerate}
\end{lemma}

\begin{proof}

(1) By Proposition \ref{prop:L1-lowerbnd},
\begin{align*}
    \int_{\R^2} Q^{s+1}\, dx  = \frac{1}{L^{(1)}_{\gamma,2}} \le \frac{4\pi \gamma^{1+\gamma}\, e^{1+2\gamma(1+\gamma_E)}}{(1+\gamma)^{1+\gamma}(1+2\gamma)} \le 4\pi e\, \gamma\, e^{2\gamma(1+\gamma_E)},
\end{align*}
since $\gamma^\gamma \le 1$ and $(1+\gamma)^{1+\gamma}(1+2\gamma)\ge 1$ for $0<\gamma\le1$. Inserting $\gamma = \frac{2}{s-1}\le \frac{3}{s}$ for $s\ge3$ shows the claim. 

(2) Follows from (1) and \eqref{eq:pohozaev}. 

(3) Follows from (2) and \eqref{eq:pohozaev3}. 

(4) Since $Q$ is radially decreasing, $\pi |x|^2 Q(x)^2 \le \int_{|y|\le |x|} Q(y)^2\, dy \le C_1/s^2$ by (2) which shows $Q(x) \le C / (s|x|)$ for all $x\ne0$.

(5) By (4), $Q(y)\le C/s$ on $\{|y|\ge 1\}$, hence $Q^s = Q^{s-1}Q \le (C/s)^{s-1}Q$ there and by \eqref{eq:pohozaev2}
\begin{align*}
    \int_{|y|\ge1} Q^s\, dy \le \Big(\frac{C}{s}\Big)^{s-1} \int_{\R^2} Q\, dy = \Big(\frac{C}{s}\Big)^{s-1} \int_{\R^2} Q^{s}\, dy .
\end{align*}
Thus, the mass of $Q^s$ concentrates in the unit disk for $s$ large. We may choose $s_0$ so large that for any $s \geq s_0$
\begin{align*}
     \int_{|y|\ge1} Q^s\, dy \le \frac{1}{2} \int_{\R^2} Q^{s}\, dy.
\end{align*}
For $|x|\ge 2$ and $|y|\le 1$ we have $|x-y|\le |x|+1$ and therefore $K_0(|x-y|)\ge K_0(|x|+1)$ by monotonicity of $K_0$. The Green representation \eqref{eq:green} gives, for $|x|\ge2$ and $s \geq s_0$,
\begin{align*}
    Q(x) \ge \frac{1}{2\pi} \int_{|y|\le1} K_0(|x-y|)\, Q(y)^s\, dy \ge \frac{\int_{\R^2} Q^{s}\, dx }{4\pi}\, K_0(|x|+1) ,
\end{align*}
hence 
\begin{align*}
    \int_{\R^2}Q^2\,dx \ge  \bigg(\frac{1}{4\pi} \int_{\R^2} Q^{s}\, dx \bigg)^2\int_{|x|\ge2}K_0(|x|+1)^2\,dx = C'\, \bigg( \int_{\R^2} Q^{s}\, dx \bigg)^2
\end{align*}
and thus $\int_{\R^2} Q^{s}\, dx\le C/s$ follows from (2).

(6) First, we just showed in (5) that $\int_{\R^2} Q\, dx = \int_{\R^2} Q^{s}\, dx \le C_1$ for large $s$. Second, by Cauchy-Schwarz and (2),
\begin{align*}
    \int_{|x|\le1} |\ln|x||\, Q(x)\, dx \le \Vert Q\Vert_{L^2(\R^2)}\, \Vert \ln|\cdot|\, \Vert_{L^2(\{|x|\le1\})} \le C_2 .
\end{align*}
Third, for the first moment, we use \eqref{eq:green} and $\int_{\R^2}|x|K_0(|x-y|)\,dx = \int_{\R^2}|x+y|K_0(|x|)\,dx \le 2\pi|y| + C_3$ to get
\begin{align*}
    \int_{\R^2} |x|\, Q(x)\, dx \le C_4 \int_{\R^2} (1+|y|)\, Q(y)^s\, dy \le 2C_4\int_{\R^2} Q^{s}\, dx + C_4\Big(\frac{C}{s}\Big)^{s-1} \int_{\R^2} |y|\, Q(y)\, dy,
\end{align*}
where in the last step we split $\int |y| Q^s\,dy$ at $|y|=1$ and used $Q^{s-1}\le (C/s)^{s-1}$ on $\{|y|\ge1\}$, as in (5). Since $\int|y|Q\,dy<\infty$ by \eqref{eq:decay}, this gives $\int_{\R^2}|x|Q\,dx\le 3C_4 \int_{\R^2} Q^{s}\, dx\le C_5$ for $s$ large. Finally, $\int_{|x|\ge1}|\ln|x||\,Q\,dx \le \int_{\R^2}|x|\,Q\,dx$. Assembling all estimates gives 
\begin{align*}
    \int_{\R^2}\big(1+|x|+|\ln|x||\big)\,Q(x)\,dx \le C.
\end{align*}

(7) The inequality follows from (5) and the identity $\int_{\R^2} Q^{s}\, dx = \frac{Q(0)}{s}I(s)$, plus the fact that $Q(0)\ge1$. 
\end{proof}

The next lemma shows that $U_s$ indeed converges locally uniformly to the limit profile $U$. We introduce as an auxiliary variable, for $r>0$, the mass of $(1+U_s/s)^s$ in the ball of radius $r$,
\begin{align} \label{eq:Irs}
    I_r(s) := \int_{|z|\le r} \Big(1+\frac{U_s}{s}\Big)^{s} \, dz ,
\end{align}
so that $I_r(s)\le I(s)$ for every $r$ and $I_r(s) \to I(s)$ for $r\to\infty$.

\begin{lemma} \label{lem:compact}
The rescaled profiles $U_s$ have the following properties:
\begin{enumerate}
\item For every $s$ and every $z\in\R^2$, one has $|U_s(z)|\le |z|^2/4$.
\item As $s\to\infty$, $U_s \to U$ and $\big(1+\frac{U_s}{s}\big)^s \to e^{U}$ locally uniformly on $\R^2$. In particular,
$$I_r(s) \to \dfrac{8\pi r^2}{8+r^2}, \qquad s\to \infty,$$
for every $r>0$.
\end{enumerate}
\end{lemma}

\begin{proof}
We write $u_s(r)$ for the radial profile of $U_s$, i.e.\ $U_s(z)=u_s(|z|)$, and
\begin{align*}
    f_s := \Big(1+\frac{u_s}{s}\Big)^{s} - s\eps^2\Big(1+\frac{u_s}{s}\Big)
\end{align*}
for the right-hand side of \eqref{eq:innerEL}. In radial coordinates, \eqref{eq:innerEL} reads 
\begin{align}
   -(r u_s')' = r f_s. \label{eq:radialform} 
\end{align}
Since $u_s(0)=0$ and $u_s'(0)=0$, integrating twice gives
\begin{align} \label{eq:compact-int}
    u_s'(r) = -\frac{1}{r}\int_0^r f_s(\tau)\,\tau\,d\tau, \qquad
    u_s(r) = -\int_0^r \frac{1}{t}\int_0^t f_s(\tau)\,\tau\,d\tau\,dt , \qquad r>0 .
\end{align}

(1) The first term of $f_s$ lies in $[0,1]$ by \eqref{eq:Qs-rescaled}, and so does the second, since $0\le 1+U_s/s = Q(\eps z)/Q(0) \le 1$ and $s\eps^2 = Q(0)^{1-s}\le 1$. Hence $|f_s|\le 1$, and \eqref{eq:compact-int} gives $|u_s'(r)|\le r/2$ and $|u_s(r)|\le r^2/4$.

(2) Let $s_n\to\infty$ be any sequence, let $\eps_n$ denote $\eps$ at $s=s_n$, and write $u_n:=u_{s_n}$, $f_n:=f_{s_n}$. By the bounds from (1), the functions $u_n$ are uniformly bounded and equicontinuous on $[0,R]$ for every $R>0$. By the Arzel\`a-Ascoli theorem and a diagonal argument, we may pass to a subsequence such that $s_n\eps_n^2\to\eta_*\in[0,1]$ and $u_n\to u_*$ locally uniformly on $[0,\infty)$. The limit $u_*$ is continuous and satisfies $u_*\le 0$ and $u_*(0)=0$.

We first show that $f_n\to f_*:=e^{u_*}-\eta_*$ locally uniformly. Since $\ln(1+x)\ge x-x^2$ for $-\frac12\le x\le 0$, we have the elementary inequality
\begin{align} \label{eq:compact-exp}
    0\le e^{a} - \Big(1+\frac{a}{s}\Big)^{s} \le \frac{a^2}{s}, \qquad -\frac{s}{2}\le a\le 0 .
\end{align}
Fix $R>0$. On $[0,R]$ we have $|u_n|\le R^2/4$ by (1), so for $n$ large \eqref{eq:compact-exp} and $|e^{a}-e^{b}|\le|a-b|$ for $a,b\le0$ give
\begin{align*}
    \Big| \Big(1+\frac{u_n}{s_n}\Big)^{s_n} - e^{u_*} \Big| \leq  \Big| \Big(1+\frac{u_n}{s_n}\Big)^{s_n} - e^{u_n} \Big| + |e^{u_n} - e^{u_*}| \le \frac{1}{s_n} \bigg(\frac{R^2}{4 } \bigg)^2 + |u_n-u_*| .
\end{align*}
Moreover, $s_n\eps_n^2\,(1+u_n/s_n)\to\eta_*$ uniformly on $[0,R]$. Hence $f_n\to f_*$ uniformly on $[0,R]$.

Passing to the limit in the second identity of \eqref{eq:compact-int}, we obtain
\begin{align} \label{eq:compact-limit}
    u_*(r) = -\int_0^r \frac{1}{t}\int_0^t \big(e^{u_*(\tau)}-\eta_*\big)\,\tau\,d\tau\,dt , \qquad r\ge 0 .
\end{align}
In particular, $u_*$ is continuously differentiable and
\begin{align*}
    r u_*'(r) = -\int_0^r e^{u_*(\tau)}\,\tau\,d\tau + \eta_*\frac{r^2}{2} .
\end{align*}

Next we show that $\eta_*=0$. Since $(1+u_n/s_n)^{s_n}\to e^{u_*}$ pointwise, Fatou's lemma and the bound $I(s)\le C$ from Lemma \ref{lem:crude} (7) give
\begin{align*}
    2\pi\int_0^\infty e^{u_*(\tau)}\,\tau\,d\tau \le \liminf_{n\to\infty} I(s_n) \le C .
\end{align*}
Suppose that $\eta_*>0$. Then
\begin{align*}
    ru_*'(r) \ge -\frac{C}{2\pi} + \eta_*\frac{r^2}{2} ,
\end{align*}
so $u_*'(r)\ge \eta_* r/4$ for all large $r$ and thus $u_*(r)\to+\infty$ as $r\to\infty$, contradicting $u_*\le0$. Hence $\eta_*=0$.

With $\eta_*=0$, equation \eqref{eq:compact-limit} states that $u_*$ is a continuous solution of
\begin{align} \label{eq:compact-fixedpoint}
    u(r) = -\int_0^r \frac{1}{t}\int_0^t e^{u(\tau)}\,\tau\,d\tau\,dt , \qquad r\ge 0 .
\end{align}
This equation has at most one continuous solution. Indeed, let $u$ and $v$ be two continuous solutions, fix $R>0$ and set $M:=\sup_{[0,R]}|u-v|$. Both solutions are nonpositive, so $|e^{u}-e^{v}|\le|u-v|$ and
\begin{align*}
    |u(r)-v(r)| \le \int_0^r \frac{1}{t}\int_0^t |u(\tau)-v(\tau)|\,\tau\,d\tau\,dt , \qquad 0\le r\le R .
\end{align*}
Inserting the bound $|u-v|\le M$ and iterating, we obtain by induction
\begin{align*}
    |u(r)-v(r)| \le \frac{M}{(k!)^2}\Big(\frac{r^2}{4}\Big)^{k} , \qquad 0\le r\le R,\ k\in\N .
\end{align*}
Letting $k\to\infty$ shows $u=v$ on $[0,R]$. The radial profile of $U$ solves \eqref{eq:compact-fixedpoint}, since $-\Delta U=e^{U}$, $U(0)=0$ and $\nabla U(0)=0$. Therefore, uniqueness implies $u_*(r)=-2\ln(1+r^2/8)$.

We have shown that every sequence $s_n\to\infty$ has a subsequence along which $U_{s_n}\to U$ locally uniformly. Hence $U_s\to U$ locally uniformly as $s\to\infty$. The locally uniform convergence of $(1+U_s/s)^s$ to $e^{U}$ follows from this and \eqref{eq:compact-exp} as above. The convergence of $I_r(s)$ follows by integrating over the bounded set $\{|z|\le r\}$ and using Lemma \ref{lem:U_integrals}.
\end{proof}

The following lemma shows that $Q(0)$ is bounded from above and stays strictly above one. This leads to exponential decay of $\varepsilon$.

\begin{lemma}\label{lem:height}
The maximum $Q(0)$ of $Q$ satisfies the following bounds:
\begin{enumerate}
\item $\limsup_{s\to\infty} Q(0) \le \sqrt{e}$.
\item There exists a constant $q_0$ such that $\liminf_{s\to\infty} Q(0) \ge q_0 > 1$.
\end{enumerate}
Consequently, there are $c>0$ and $s_0$ such that for $s\ge s_0$
\begin{align} \label{eq:eps-bounds}
    e^{-s} \le \eps \le e^{-cs}, \qquad s\eps^2 \le e^{-cs} .
\end{align}
\end{lemma}

\begin{proof}
(1) Fix $r>0$. Since $K_0$ is positive and decreasing, \eqref{eq:I-exact} gives
\begin{align*}
    2\pi s \ge \int_{|z|\le r} K_0(\eps|z|)\Big(1+\frac{U_s}{s}\Big)^s\, dz \ge K_0(\eps r)\, I_r(s) .
\end{align*}
Since $Q(0)\geq 1$, $\eps\le s^{-1/2}\to0$, so \eqref{eq:K0-expansion} gives $K_0(\eps r) = \ln\frac1\eps + \ln\frac2r - \gamma_E + o(1)$, and $I_r(s)\to\frac{8\pi r^2}{8+r^2}$ for $s\to \infty$ by Lemma \ref{lem:compact} (2). Therefore
\begin{align*}
    \ln\frac1\eps \le \frac{2\pi s}{I_r(s)} + \ln\frac{r}{2} + \gamma_E + o(1) = \frac{s(8+r^2)}{4 r^2} (1+o(1) )  , \qquad s\to \infty,
\end{align*}
and by the definition \eqref{eq:def_rescaled_var} of $\eps$,
\begin{align*}
    (s-1)\ln Q(0) = 2\ln\frac1\eps - \ln s \le \frac{s}{2}\Big(1+\frac{8}{r^2}\Big)(1+o(1)) .
\end{align*}
Hence $\limsup_{s\to\infty}\ln Q(0) \le \frac12 (1+8/r^2)$ for every $r>0$, which proves (1).

(2) By Lemma \ref{lem:K0}, $K_0(\rho)\le \ln^+\frac{2}{\rho} + C$ for all $\rho>0$. Inserting this into \eqref{eq:I-exact} and using $(1+U_s/s)^s\le1$, $\ln\frac2\eps>0$ and $I(s)\le C_0$,
\begin{align*}
    2\pi s &\le \int_{\eps|z|\le 2} \Big(\ln\frac{2}{\eps} - \ln|z|\Big)\Big(1+\frac{U_s}{s}\Big)^s\, dz + C I(s) \\
    &\le \ln\frac2\eps\, I(s) + \int_{|z|\le1} \ln\frac{1}{|z|}\, dz + C I(s) \le C_0 \ln\frac2\eps + C .
\end{align*}
Thus $\ln\frac1\eps \ge \frac{2\pi s}{C_0} - C$, and \eqref{eq:def_rescaled_var} gives $(s-1)\ln Q(0) = 2\ln\frac1\eps - \ln s \ge \frac{4\pi s}{C_0} - \ln s - C$, which proves (2).

For \eqref{eq:eps-bounds}, we use $2\ln\frac1\eps = (s-1)\ln Q(0) + \ln s$. By (1), $\ln\frac1\eps\le s$ for large $s$, while by (2), $\ln Q(0)\ge 2\pi/C_0$ for large $s$, so $s\eps^2 = Q(0)^{1-s}\le e^{-2\pi(s-1)/C_0}$.
\end{proof}

We now come to the decay estimate. Fix $\theta\in(0,\frac12)$ and set
\begin{align} \label{eq:rs}
    r_s := \bigg(\frac{\theta}{s\eps^2}\bigg)^{1/2} = \frac{1}{\eps}\bigg(\frac{\theta}{s}\bigg)^{1/2}, \qquad \delta_s := \eps r_s = \bigg(\frac{\theta}{s}\bigg)^{1/2} .
\end{align}
In the original variable $x=\eps z$, the ball $\{|z|\le r_s\}$ corresponds to $\{|x|\le\delta_s\}$, and $r_s\ge e^{cs}$ by \eqref{eq:eps-bounds}. On $\{|z|\le r_s\}$ the forcing term in \eqref{eq:innerEL} is negligible, and $U_s$ is controlled by comparison with the Liouville profile. Outside, on $\{|z|\ge r_s\}$, $Q$ is already small and $Q^s$ is superexponentially small.

\begin{lemma}  \label{lem:decay}
There exist $R=R_\theta \geq 1$, $s_\theta$ and $C=C_\theta$ independent of $s$ such that for all $s\ge s_\theta$:
\begin{enumerate}
\item $-C (1+\ln|z| )\le U_s(z) \le 0$ for $|z|\ge1$,
\item $\Big(1+\dfrac{U_s(z)}{s}\Big)^s\le \Big(\dfrac{R}{|z|}\Big)^{4-2\theta}$ for $R\le|z|\le r_s$,
\item $Q(x)\le C(\theta s)^{-1/2}$ for $|x|\ge\delta_s$, and
\begin{align*}
    \int_{|z|\ge r_s} (1+|z|+|\ln|z||+|U_s(z)| )\Big(1+\frac{U_s(z)}{s}\Big)^s\,dz\ \le\ s^{-s/4} .
\end{align*}
\end{enumerate}
\end{lemma}

\begin{proof}
(1) We only need to show the lower bound. Integrating the radial form \eqref{eq:radialform} from $0$ to $r$, we obtain the identity
\begin{align} \label{eq:radial}
    r u_s'(r) = -\frac{I_r(s)}{2\pi} + s\eps^2\int_0^r\Big(1+\frac{u_s(t)}{s}\Big)\,t\,dt .
\end{align}
Dropping the nonnegative second term in \eqref{eq:radial} and using $I_r(s)\le I(s)\le C$ from Lemma \ref{lem:crude} (7) gives $ru_s'(r)\ge -\frac{C}{2\pi}$, hence $u_s(r)\ge u_s(1) - \frac{C}{2\pi}\ln r \ge -\frac14 - \frac{C}{2\pi}\ln r$ for $r\ge1$, by Lemma \ref{lem:compact} (1).

(2) By Lemma \ref{lem:compact} (2) and Lemma \ref{lem:U_integrals}, we can choose $R$ so large that $\frac{8\pi R^2}{8+R^2}\ge 2\pi(4-\frac{\theta}{2})$, and then $s_\theta$ so large that $I_R(s)\ge 2\pi(4-\theta)$ for all $s\ge s_\theta$. Since $I_r(s)$ is nondecreasing in $r$, $I_r(s)\ge2\pi(4-\theta)$ for $r\ge R$. On the other hand, $0\le 1+U_s/s\le 1$ gives for the second term in \eqref{eq:radial} the bound 
\begin{align*}
     s\eps^2\int_0^r\Big(1+\frac{u_s(t)}{s}\Big)\,t\,dt \leq s\eps^2 \frac{ r^2}{2}\le s\eps^2 \frac{r_s^2}{2} = \frac{ \theta}{2}, \qquad r\le r_s,
\end{align*}
by the definition \eqref{eq:rs} of $r_s$. Hence, for $R\le r\le r_s$,
\begin{align*}
    ru_s'(r) \le -(4-\theta) + \frac\theta2 \le -(4-2\theta) .
\end{align*}
Integrating from $R$ to $r$ and using $u_s(R)\le0$ gives $u_s(r)\le -(4-2\theta)\ln\frac rR$, and the bound on $(1+U_s/s)^s$ follows from \eqref{eq:Qs-rescaled}.

(3) By Lemma \ref{lem:crude} (4), $Q(x)\le \frac{C}{s\delta_s} = C(\theta s)^{-1/2}$ for $|x|\ge\delta_s$, and hence
\begin{align*}
    Q(x)^s\le \big(C(\theta s)^{-1/2}\big)^{s-1}\,Q(x), \qquad |x|\ge \delta_s .
\end{align*}
Returning to the variable $z$ via $x=\eps z$, \eqref{eq:def_rescaled_var}, $Q(\eps z)^s = \frac{Q(0)}{s\eps^2}(1+U_s(z)/s)^s$ and using $Q(0)>1$, we obtain for any weight $w\ge0$
\begin{align*}
    \int_{|z|\ge r_s} w(z) \Big(1+\frac{U_s(z)}{s}\Big)^s\, dz
    &= \frac{s}{Q(0)}\int_{|x|\ge\delta_s} w\Big(\frac{x}{\eps}\Big)\, Q(x)^s\, dx \\
    &\le s (C(\theta s)^{-1/2})^{s-1} \int_{\R^2} w\Big(\frac{x}{\eps}\Big)\, Q(x)\, dx .
\end{align*}
For $w(z)=1$, $|z|$, $|\ln|z||$ and $|U_s(z)|$ we have, respectively, $w(x/\eps) \le 1$, $|x|/\eps$, $|\ln|x||+\ln\frac1\eps$ and $s$, so that by Lemma \ref{lem:crude} (6) and $\ln\frac1\eps\le s$, $\frac1\eps\le e^{s}$ from \eqref{eq:eps-bounds}, the last integral is bounded by $Ce^{s}$ for large $s$. Altogether the tail integral is at most
\begin{align*}
    C s\,(s+e^{s})\,\big(C(\theta s)^{-1/2}\big)^{s-1} = \exp\Big(-\frac{s-1}{2}\ln s + O(s)\Big) \le s^{-s/4}
\end{align*}
for $s\ge s_\theta$.
\end{proof}

\begin{proof}[Proof of Proposition \ref{prop:IJL_limits}]
Let $\theta\in(0,\frac12)$ and $R$, $s_\theta$, $r_s$ be as in Lemma \ref{lem:decay}. We first build an $s$-independent majorant on $\{|z|\le r_s\}$. By \eqref{eq:Qs-rescaled} and Lemma \ref{lem:decay} (2),
\begin{align*}
    \Big(1+\frac{U_s(z)}{s}\Big)^s \le g_0(z) := \min\Big\{1, \Big(\frac{R}{|z|}\Big)^{4-2\theta}\Big\}, \qquad |z|\le r_s,\ s\ge s_\theta ,
\end{align*}
and by Lemma \ref{lem:compact} (1) and Lemma \ref{lem:decay} (1), 
\begin{align*}
|U_s(z)|\le C(1+\ln(1+|z|)), \qquad z\in\R^2,\ s\ge s_\theta.
\end{align*}
Hence the three integrands $(1+U_s/s)^s$, $U_s(1+U_s/s)^s$ and $\ln|z|\,(1+U_s/s)^s$, restricted to $\{|z|\le r_s\}$, are bounded in absolute value by
\begin{align*}
    g(z) := g_0(z)\,\big(1 + |\ln|z|| + C(1+\ln(1+|z|))\big),
\end{align*}
which belongs to $L^1(\R^2)$ because $|\ln|z||$ is integrable near the origin and $4-2\theta>2$. By Lemma \ref{lem:compact} (2), the three integrands of $I$, $J$ and $L$ converge pointwise to $e^{U}$, $Ue^{U}$ and $\ln|z|\,e^{U}$, respectively, and the indicator functions of $\{|z|\le r_s\}$ converge pointwise to $1$ since $r_s\to\infty$. Dominated convergence therefore gives
\begin{align*}
    \int_{|z|\le r_s}\Big(1+\frac{U_s}{s}\Big)^s dz&\to\int_{\R^2}e^{U}dz, &
    \int_{|z|\le r_s}U_s\Big(1+\frac{U_s}{s}\Big)^s dz&\to\int_{\R^2}Ue^{U}dz, \\
    \int_{|z|\le r_s}\ln|z|\Big(1+\frac{U_s}{s}\Big)^s dz&\to\int_{\R^2}\ln|z|\,e^{U}dz .
\end{align*}
The contributions of $\{|z|\ge r_s\}$ to $I(s)$, $J(s)$ and $L(s)$ are bounded by $s^{-s/4}$ by Lemma \ref{lem:decay} (3). Together with the values of the limiting integrals from Lemma \ref{lem:U_integrals}, this proves Proposition \ref{prop:IJL_limits}.
\end{proof}

\begin{proof}[Proof of Lemma \ref{lem:remainders}]
We keep $\theta$, $R$, $r_s$, $\delta_s$ from Lemma \ref{lem:decay} and consider $s\ge s_\theta$.

\emph{The remainder $R_1(s)$.} On $\{|z|\le r_s\}$ we have $\eps|z|\le\delta_s\le1$, so by \eqref{eq:K0-expansion} the integrand in \eqref{eq:R1-def} is nonnegative and bounded by
\begin{align*}
    C\eps^2|z|^2  \Big(1+\ln\frac{1}{\eps|z|} \Big)\,\Big(1+\frac{U_s}{s}\Big)^s
    &\le C\eps^2|z|^2\Big(1+\ln\frac1\eps+\ln^+\frac1{|z|}\Big)\Big(1+\frac{U_s}{s}\Big)^s \\
    &\le C' s \eps^2 |z|^2 \Big(1+\ln^+\frac1{|z|}\Big) g_0(z),
\end{align*}
where we used \eqref{eq:eps-bounds} and the majorant $g_0$ from the previous proof. Hence
\begin{align*}
    & \int_{|z|\le r_s}\Big[K_0(\eps|z|) - \ln\frac{2}{\eps|z|}+\gamma_E\Big]\Big(1+\frac{U_s}{s}\Big)^s dz\\
    & \qquad \leq \bigg( \int_{|z|\le R} + \int_{R\leq |z|\le r_s} \bigg) \, C s\eps^2 |z|^2 \Big(1+\ln^+\frac1{|z|}\Big) g_0(z) dz\\
    &\qquad \le C s\eps^2\bigg(R^4 + R^{4-2\theta}\int_R^{r_s} t^{-1+2\theta}\,dt\bigg) \\
    &\qquad \le C' s\eps^2 r_s^{2\theta} \le C s \eps^{2-2\theta},
\end{align*}
since $r_s^{2\theta} = \theta^\theta s^{-\theta}\eps^{-2\theta}\le \eps^{-2\theta}$. This is $O(e^{-cs})$ by \eqref{eq:eps-bounds}. On $\{|z|\ge r_s\}$ we bound the integrand in \eqref{eq:R1-def} by $\big(K_0(\eps|z|) + \ln\frac1\eps + |\ln|z|| + C\big)(1+U_s/s)^s$. The terms with $\ln\frac1\eps\le s$, $|\ln|z||$ and $C$ contribute at most $C s\cdot s^{-s/4}$ by Lemma \ref{lem:decay} (3). For the term with $K_0$ we return to the variable $x=\eps z$ as in the proof of Lemma \ref{lem:decay} (3). This gives
\begin{align*}
    \int_{|z|\ge r_s} K_0(\eps|z|)\Big(1+\frac{U_s}{s}\Big)^s dz &= \frac{s}{Q(0)}\int_{|x|\ge\delta_s}K_0(|x|)\,Q(x)^s\,dx \\
    &\le s (C(\theta s)^{-1/2} )^{s-1}\, Q(0)\int_{\R^2}K_0(|x|)\,dx ,
\end{align*}
which is $O(s^{-s/4})$ since $\int_{\R^2}K_0(|x|)\,dx = 2\pi$ and $Q(0)\le C$ by Lemma \ref{lem:height} (1). Altogether, $R_1(s) = O(e^{-cs})$.

\emph{The remainder $R_2(s)$.} We first show that there exists a constant $C$ such that
\begin{align} \label{eq:Ds-bound}
    R_2(s) \le \frac{C}{\int_{\R^2} Q^{s} \, dx } \int_{\R^2} m(y) \, Q(y)^s \, dy, \quad m(y):=\min\Big\{ |y| \Big(1+\ln^{+}\frac{1}{|y|} \Big)^{1/2},\, 1\Big\} .
\end{align}
By \eqref{eq:green},
\begin{align*}
    Q(x) - \frac{\int_{\R^2} Q^s \, dy}{2\pi} K_0(|x|) = \frac{1}{2\pi} \int_{\R^2} \big[K_0(|x-y|) - K_0(|x|)\big] Q(y)^s \, dy .
\end{align*}
By Minkowski's integral inequality, the $L^2$ norm of the left-hand side is at most
\begin{align*}
    \frac{1}{2\pi} \int_{\R^2} \Vert K_0(|\cdot - y|) - K_0(|\cdot|)\Vert_{L^2(\R^2)}\, Q(y)^s \, dy ,
\end{align*}
so \eqref{eq:Ds-bound} follows from the definition \eqref{eq:Ds-def} of $R_2(s)$ once we show
\begin{align} \label{eq:K0-translation}
    \Vert K_0(|\cdot - y|) - K_0(|\cdot|)\Vert_{L^2(\R^2)} \le C\, m(y), \qquad y\in\R^2 .
\end{align}
Since $\frac{1}{2\pi}K_0(|\cdot|)$ is the fundamental solution of $-\Delta+1$, Plancherel's theorem gives
\begin{align*}
    \Vert K_0(|\cdot - y|) - K_0(|\cdot|)\Vert_{L^2(\R^2)}^2 = \int_{\R^2} \frac{|e^{-iy\cdot\xi}-1|^2}{(1+|\xi|^2)^2}\, d\xi .
\end{align*}
We use $|e^{-iy\cdot\xi}-1|^2\le\min\{|y|^2|\xi|^2,4\}$. The constant upper bound by $4$ shows that the right-hand side is at most $4\pi$ for every $y$, so \eqref{eq:K0-translation} holds in particular for $|y|\geq 1$. For $|y|\le1$, we split the integral at $|\xi|=1/|y|$. Then
\begin{align*}
    \Vert K_0(|\cdot - y|) - K_0(|\cdot|)\Vert_{L^2(\R^2)}^2 &\leq |y|^2\int_{|\xi|\le 1/|y|} \frac{|\xi|^2}{(1+|\xi|^2)^2}\, d\xi + \int_{|\xi|> 1/|y|} \frac{4}{(1+|\xi|^2)^2}\, d\xi 
\end{align*}
where
\begin{align*}
 |y|^2\int_{|\xi|\le 1/|y|} \frac{|\xi|^2}{(1+|\xi|^2)^2}\, d\xi   
    = 2\pi |y|^2 \int_0^{1/|y|} \frac{r^3}{(1+r^2)^2}\, dr & = \pi |y|^2 \int_0^{1/|y|^2} \frac{u}{(1+u)^2}\, du \\
    &\le \pi|y|^2 \ln\Big(1+\frac{1}{|y|^2}\Big) \\
    &\le C|y|^2\Big(1+\ln\frac{1}{|y|}\Big)
\end{align*}
and
\begin{align*}   
    \int_{|\xi|> 1/|y|} \frac{4}{(1+|\xi|^2)^2}\, d\xi &\le 8\pi\int_{1/|y|}^\infty r^{-3}\, dr = 4\pi |y|^2 .
\end{align*}
These bounds prove \eqref{eq:K0-translation} also for $|y|\le1$.
 
Now split the integral on the right-hand side of \eqref{eq:Ds-bound} at $|y|=\delta_s$:
\begin{align*} 
    R_2(s) \le \frac{C}{\int_{\R^2} Q^{s} \, dx } \bigg(\int_{|y|\le\delta_s}+ \int_{|y|\ge\delta_s} \bigg) \, m(y) \, Q(y)^s \, dy.
\end{align*}
For $|y|\ge\delta_s$ we use $m(y)\leq 1$, Lemma \ref{lem:decay} (3) and \eqref{eq:pohozaev2}, giving
\begin{align*}
    \int_{|y|\ge\delta_s}Q(y)^s\,dy \le \big(C(\theta s)^{-1/2}\big)^{s-1}\int_{\R^2}Q\,dy = \big(C(\theta s)^{-1/2}\big)^{s-1} \int_{\R^2} Q^s \, dy .
\end{align*}
For $|y|\le\delta_s$ we substitute $y=\eps z$, $|z|\le r_s$. Since $1+\ln\frac{1}{|y|} = 1+\ln\frac1\eps+\ln\frac1{|z|}\le C s\,(1+\ln^+\frac1{|z|})$ by \eqref{eq:eps-bounds}, and $|z|(1+\ln\frac1{|z|})^{1/2}\le C$ for $|z|\le 1$, we have $|y|(1+\ln\frac{1}{|y|})^{1/2}\le C\sqrt s\eps\,(1+|z|)$ for $|z|\le r_s$. Therefore, using $Q(y)^s\,dy = \frac{Q(0)}{s}(1+U_s/s)^s\,dz$,
\begin{align*}
    \int_{|y|\le\delta_s}|y| \Big(1+\ln^{+}\frac{1}{|y|} \Big)^{1/2}\,Q(y)^s\,dy
    &\le C\sqrt s\,\eps\,\frac{Q(0)}{s}\int_{|z|\le r_s}(1+|z|)\Big(1+\frac{U_s}{s}\Big)^s dz \\
    &\le C'\sqrt s\,\eps\,\frac{Q(0)}{s} ,
\end{align*}
because $\int_{|z|\le r_s}(1+|z|)(1+U_s/s)^s\,dz\le\int_{\R^2}(1+|z|)g_0(z)\,dz<\infty$, as $4-2\theta>3$. Finally, $\frac{Q(0)}{s} = \frac{1}{I(s)}\int_{\R^2} Q^s \, dy$ and $I(s) \to 8 \pi$ by Proposition \ref{prop:IJL_limits}, so that $\frac{Q(0)}{s}\le C \int_{\R^2} Q^s \, dy$ for large $s$. Inserting both parts into \eqref{eq:Ds-bound} gives
\begin{align*}
    R_2(s) \le C\Big(\sqrt{s}\,\eps + \big(C(\theta s)^{-1/2}\big)^{s-1}\Big) = O(e^{-cs}), \qquad s \to \infty,
\end{align*}
by \eqref{eq:eps-bounds}. This completes the proof of Lemma \ref{lem:remainders}.
\end{proof}

\subsection{Proof of Theorem \ref{thm:L1thm}}

With Proposition \ref{prop:IJL_limits} and Lemma \ref{lem:remainders}, we can now prove the main theorem. We first prove a second order term for $I$ using Lemma \ref{lem:P}.

\begin{lemma} \label{lem:I_limit_impr}
   \begin{align}
            I(s)=8\pi\Big(1-\frac3s+o\Big(\frac1s\Big)\Big), \qquad s\to\infty.
\end{align} 
\end{lemma}

\begin{proof}
Because of Proposition \ref{prop:IJL_limits}, we can write $ I(s)=8\pi (1+\eta(s))$ with $\eta(s) \to 0$ as $s\to \infty$. Now the relation \eqref{eq:P-inner} from Lemma \ref{lem:P} together with the limits \eqref{eq:I_limit} and \eqref{eq:J_limit} for $I$ and $J$ and the exponential smallness of $R_2(s)$ from Lemma \ref{lem:remainders} implies
\begin{align*}
    8\pi (1+\eta(s)) - \frac{16\pi}{s} + o\bigg(\frac{1}{s} \bigg) = 8\pi \bigg(1 +\frac{1}{s} \bigg) (1+\eta(s) )^2 .
\end{align*}
Thus,
\begin{align*}
    -8\pi \eta(s) (1+o(1)) = \frac{24\pi}{s} (1+o(1))
\end{align*}
and therefore 
\begin{align*}
    \eta(s) = - \frac{3}{s} + o \bigg(\frac{1}{s} \bigg).
\end{align*}
\end{proof}

From Lemma \ref{lem:I}, we can now deduce an asymptotic expansion for $\varepsilon$.

\begin{lemma} \label{lem:scalebubble}
    \begin{align*}
        \ln\frac{1}{\eps} = \frac{s}{4} + \frac34 + \frac12\ln 2 + \gamma_E + o(1), \qquad s\to\infty .
    \end{align*}
\end{lemma}

\begin{proof}
    By Lemma \ref{lem:I}, 
    \begin{align*}
        \ln \frac{1}{\varepsilon} = \gamma_E - \ln 2 + \frac{2\pi s + L(s) - R_1(s)}{I(s)} .
    \end{align*}
Using the improved asymptotic expansion for $I$ from Lemma \ref{lem:I_limit_impr}, the limit \eqref{eq:L_limit} for $L$ from Proposition \ref{prop:IJL_limits} and the exponential smallness of $R_1(s)$ from Lemma \ref{lem:remainders}, we obtain
    \begin{align*}
        \ln \frac{1}{\varepsilon} &= \gamma_E - \ln 2 + \frac{2\pi s + 12 \pi \ln 2 + o(1)}{8\pi} \bigg(1 + \frac{3}{s} +o \bigg(\frac{1}{s} \bigg) \bigg) \\
        &= \frac{s}{4} + \frac34 + \frac12\ln 2 + \gamma_E + o(1), \qquad s\to\infty .
    \end{align*} 
\end{proof}

Finally, we can turn the asymptotic expansion for $\varepsilon$ into one for $Q(0)$ and together with the asymptotic expansion for $I$, we obtain asymptotic expansions for different integral moments of $Q$.

\begin{proposition} \label{prop:norms} As $s\to\infty$,
\begin{enumerate}
    \item  \begin{align}
        Q(0)^2  &= e \bigg( 1- 2\frac{\ln s}{s} + (4+4\gamma_E + 2\ln 2) \frac{1}{s} + o\Big(\frac{1}{s}\Big) \bigg).
    \end{align}
    \item \begin{align}
    \int_{\R^2} Q^{s+1} \, dx &= \frac{8\pi e}{s}\bigg(1 - 2\frac{\ln s}{s} + (-1 + 4\gamma_E + 2\ln 2)\frac1s + o\bigg(\frac1s\bigg)\bigg),  \label{eq:L2p-asympt}\\
    \int_{\R^2} Q^s \, dx &= \frac{8\pi\sqrt e}{s}\bigg(1 - \frac{\ln s}{s} + (-1 + 2\gamma_E + \ln 2)\frac1s + o\bigg(\frac1s\bigg)\bigg), \label{eq:Ls-asympt}\\
    \int_{\R^2} Q^2 \, dx &= \frac{16\pi e}{s^2}\bigg(1 - 2\frac{\ln s}{s} + (-2 + 4\gamma_E + 2\ln 2)\frac1s + o\bigg(\frac1s\bigg)\bigg).  \label{eq:L2-asympt}
\end{align}
\end{enumerate}
\end{proposition}

\begin{proof}
    (1) By definition of $\varepsilon$, see \eqref{eq:def_rescaled_var}, and Lemma \ref{lem:scalebubble},
    \begin{align*}
        (s-1) \ln Q(0) = 2 \ln \frac{1}{\varepsilon} - \ln s = \frac{s}{2} - \ln s + \frac{3}{2} + \ln 2 + 2\gamma_E  + o(1), \qquad s\to \infty,
    \end{align*}
    and therefore
    \begin{align*}
     \ln Q(0) =\frac{1}{2} - \frac{\ln s}{s} + (2 + \ln 2 + 2\gamma_E)\frac{1}{s}    + o \bigg(\frac{1}{s} \bigg), \qquad s\to \infty,
    \end{align*}
    and the desired asymptotic expansion follows.

    (2) Since $\int_{\R^2}  Q^{s} dx = \frac{Q(0)}{s} I(s)$, Lemma \ref{lem:I_limit_impr} and (1) give
    \begin{align*}
        \int_{\R^2}  Q^{s} dx =  \frac{8\pi \sqrt{e} }{s} \bigg( 1- \frac{\ln s}{s} + (2+2\gamma_E + \ln 2)\frac{1}{s}  \bigg) \Big(1-\frac3s\Big)  + o\bigg(\frac{1}{s^2}\bigg), \qquad s\to \infty,
    \end{align*}
    from which \eqref{eq:Ls-asympt} follows. By \eqref{eq:P} from Lemma \ref{lem:P} and the exponential smallness of $R_2(s)$ from Lemma \ref{lem:remainders}, we then obtain \eqref{eq:L2p-asympt}. Finally, \eqref{eq:L2-asympt} follows from \eqref{eq:pohozaev} and \eqref{eq:L2p-asympt}.    
\end{proof}

\begin{proof}[Proof of Theorem \ref{thm:L1thm}]
    The asymptotics follow from Proposition \ref{prop:Q} (4) and \eqref{eq:L2p-asympt}. Indeed, 
    \begin{align*}
        L_{\gamma,2}^{(1)} = \frac{1}{ \int_{\R^2} Q^{s+1} \, dx} &= \frac{s}{8\pi e}\bigg(1 + 2\frac{\ln s}{s} - (-1 + 4\gamma_E + 2\ln 2)\frac1s + o\bigg(\frac1s\bigg)\bigg) \\
        &=\frac{1}{8\pi e}\bigg(s + 2\ln s  + 1 - 4\gamma_E - 2\ln 2 + o(1)\bigg) 
    \end{align*}
    and plugging in $s=1 +\frac{2}{\gamma}$ and using $2\ln s - 2 \ln 2 = 2 \ln (\frac{1}{\gamma} + \frac{1}{2}) = 2 \ln (\frac{1}{\gamma}) + O(\gamma)$ finishes the proof for the asymptotic expansion of $L_{\gamma,2}^{(1)}$. Proposition \ref{prop:Q} (4) also states that the corresponding optimal potential realizing $L_{\gamma,2}^{(1)}$ in \eqref{eq:L1} is, up to rescaling and symmetries, given by $Q^{s-1}$. The quotient
    \begin{align*}
        \frac{|\lambda_1|^{\gamma}}{\int_{\R^2}V^{1+\gamma}\, dx}
    \end{align*}
    is invariant under the scaling $V \to \alpha^2 V(\alpha \, \cdot \, )$, $\alpha>0$, so
    \begin{align*}
    V_\gamma(x) = 8\varepsilon^{2} Q( \sqrt{8}\varepsilon x)^{s-1}=\frac{8}{s}\bigg(1+\frac{U_s(\sqrt{8}x)}{s}\bigg)^{s-1}
    \end{align*}
    is also an optimal potential. Since $8/s = 8\gamma /(2+\gamma) = 4\gamma (1+o(1))$ as $\gamma\to 0$ and
    by Lemma \ref{lem:compact} (2),
    \begin{align*}
        \bigg(1+\frac{U_s(\sqrt{8}x)}{s}\bigg)^{s-1} \to e^{U(\sqrt{8}x)} = (1+|x|^2)^{-2}, \qquad s\to \infty,
    \end{align*}
    locally uniformly, we find
    \begin{align*}
    V_\gamma(x) = 4\gamma (1+|x|^2)^{-2} (1+o(1)) , \qquad \gamma\to 0,
    \end{align*}
    locally uniformly. This finishes the proof.    
    \end{proof}

\section*{AI usage declaration}

This paper was developed and written within 10 days, through intensive discussions and interactions with Claude Fable 5.1 and 5.2. While I guided the AI and have rewritten and reorganized most of the AI's output for the presentation in this paper, almost all of the mathematical reasoning is due to the AI, which I have checked independently while writing up the paper. But as such, I do not feel comfortable submitting this paper to a journal. Still, I want to share the result and I take full responsibility for the content of the article.

\def\myarXiv#1#2{\href{http://arxiv.org/abs/#1}{\texttt{arXiv:#1\, [#2]}}}

\end{document}